\documentclass[10pt]{amsart}
\usepackage[top=2cm,bottom=2cm,left=3cm,right=3cm,marginparwidth=2cm]{geometry}

\usepackage{mathtools} % loads amsmath
\usepackage{amssymb,amsthm}
\usepackage{xcolor}
\usepackage[square,compress,comma,numbers,sort]{natbib}

\usepackage[colorlinks=true,citecolor=blue,linkcolor=blue,urlcolor=blue]{hyperref}
\usepackage[shortlabels]{enumitem}

\numberwithin{equation}{section}

\theoremstyle{plain}
\newtheorem{theorem}{Theorem}[section]
\newtheorem{lemma}[theorem]{Lemma}

\theoremstyle{definition}

\theoremstyle{remark}
\newtheorem{remark}[theorem]{Remark}

\newcommand{\R}{\mathbb{R}}
\newcommand{\Sgen}{S^{d-1}_+(\|\cdot\|)}
\DeclarePairedDelimiterXPP\E[1]{\mathbb{E}}\{ \}{}{#1}
\DeclarePairedDelimiterXPP\pk[1]{\mathbb{P}}\{ \}{}{#1}
\DeclarePairedDelimiterXPP\ind[1]{\mathbb{I}}( ){}{#1}
\newcommand{\norm}[1]{\lVert#1\rVert}
\newcommand{\TT}{\mathcal{T}}
\newcommand{\dK}{\mathrm{d}_{K}}
\newcommand{\dW}{W_{1,\|\cdot\|_\infty}}
\newcommand{\dWd}{W}
\newcommand{\Law}{\mathrm{Law}}
\begin{document}

\title[Distances between max-stable distributions]
{Kolmogorov and Wasserstein Distances between Max-Stable Distributions}

\author{Enkelejd Hashorva}
\address{Enkelejd Hashorva, HEC Lausanne,
  University of Lausanne,\newline
  Chamberonne, 1015 Lausanne, Switzerland}
\email{Enkelejd.Hashorva@unil.ch}
\date{}

\begin{abstract}
We derive explicit comparison bounds for multivariate max-stable distributions
with unit-$\alpha$-Fr\'echet margins. For the Kolmogorov distance, the bounds
are expressed through Wasserstein distances between powered de Haan
representers, total variation distances between angular measures, and
discrepancies of the $\Psi$-functions in the inf--argmax decomposition. On the
positive $\ell_\alpha$-sphere, the coefficient multiplying the setwise
angular total-variation distance contains no explicit dimension factor for
the unnormalised angular measures used here. Separately, for
$1\le p<\alpha$, a synchronous de Haan--LePage coupling bounds the
$p$-Wasserstein distance between the max-stable laws by an
$\alpha$-Wasserstein transport cost between their unpowered de Haan
representers. We also compare laws with a common
extreme-value copula and different Fr\'echet indices, obtaining an exact
$\ell_1$-Wasserstein formula when $p=1$, and discuss applications to Archimax
and clustered Archimax copulas and to Brown--Resnick/H\"usler--Reiss models.
\end{abstract}


\maketitle

\section{Introduction}

Let \(F_1,F_2\) be two \(d\)-dimensional max-stable distribution functions.
A basic problem, relevant both for approximation theory and for statistical
modelling of extremes, is to quantify how far \(F_1\) and \(F_2\) are from one
another in terms of their spectral or parametric representations. In this
paper we study this question mainly through the Kolmogorov distance
\[
  \dK(F_1,F_2)
  :=
  \sup_{x\in\R^d}|F_1(x)-F_2(x)|,
\]
and, when finite moments are available, through Wasserstein distances between
the corresponding max-stable laws.
Recent Stein-method bounds in Kolmogorov and Wasserstein distances, allowing
the stability index and/or angular measure to vary, were obtained by
Costac\`eque and Decreusefond \cite{costaceque2025stein}.

Unless stated otherwise, we work with unit-\(\alpha\)-Fréchet margins,
\[
  F_{k,j}(x)=\exp\{-x^{-\alpha}\},
  \qquad x>0,\quad j=1,\ldots,d,\quad k=1,2,
\]
where \(\alpha>0\). If \(Z^{(k)}=(Z^{(k)}_1,\ldots,Z^{(k)}_d)\) is a de Haan
representer of \(F_k\), normalised by
\[
  \E{(Z^{(k)}_j)^\alpha}=1,
  \qquad j=1,\ldots,d,
\]
then the exponent function \(V_k=-\log F_k\) is given by
\[
  V_k(x)
  =
  \E*{
    \max_{1\le j\le d}
    \frac{(Z^{(k)}_j)^\alpha}{x_j^\alpha}
  },
  \qquad x\in(0,\infty)^d.
\]
See \cite{deHaan,Haan2006} for the spectral representation and
\cite{ressel2013homogeneous,FalkD} for the associated homogeneous
representations.
Equivalently, \(V_k\) may be represented through an angular measure \(H_k\)
on a positive unit sphere, or through the stable tail dependence function
\[
  \ell_k(z):=V_k(z^{-1/\alpha}),
  \qquad z\in[0,\infty)^d.
\]

The starting point of the paper is the observation that, after exploiting
homogeneity, the Kolmogorov distance can be reduced to a uniform comparison
of exponent functions on a homogeneous cross-section. More precisely, with
\[
  \omega(t)
  :=
  \begin{cases}
  0, & t=0,\\[1mm]
  \displaystyle
  \frac{t}{1+t}(1+t)^{-1/t}, & t>0,
  \end{cases}
\]
we obtain the intrinsic bound
\[
  \dK(F_1,F_2)
  \le
  \omega(\Delta_V)
  \le
  \frac{\Delta_V}{e},
\]
where
\[
  \Delta_V
  :=
  \sup_{\substack{x\in(0,\infty)^d\\ \min_j x_j=1}}
  |V_1(x)-V_2(x)|.
\]
This reduction is the common mechanism behind all Kolmogorov estimates in the
paper.

Our first group of results bounds \(\Delta_V\) directly in terms of de Haan
representers and angular measures. We prove, for instance,
\[
  \dK(F_1,F_2)
  \le
  \omega\left(
    \mathcal W_\alpha
    \wedge
    M_\infty\|H_1-H_2\|_{\mathrm{TV}}
  \right),
\]
where
\[
  \mathcal W_\alpha
  :=
  \inf_{Z^{(1)},Z^{(2)}}
  W_{1,\|\cdot\|_\infty}
  \bigl((Z^{(1)})^\alpha,(Z^{(2)})^\alpha\bigr),
\]
the infimum being over all de Haan representers of \(F_1\) and \(F_2\), and
\[
  M_\infty
  :=
  \sup_{s\in\Sgen}\|s\|_\infty^\alpha.
\]
In particular,
\[
  \dK(F_1,F_2)
  \le
  \frac1e
  W_{1,\|\cdot\|_\infty}
  \bigl((Z^{(1)})^\alpha,(Z^{(2)})^\alpha\bigr)
\]
for any fixed pair of representers. Since de Haan representers are not
unique, we also formulate the bound in terms of canonical representers
obtained from angular measures. On the positive \(\ell_\alpha\)-sphere this
gives \(M_\infty=1\), so, for the unnormalised angular measures used here,
the coefficient multiplying the setwise angular total-variation distance
contains no explicit dimension factor; compare
\cite{costaceque2025stein} after converting angular and total-variation
conventions.

A second set of results gives a different kind of Wasserstein control between
the max-stable laws themselves. These estimates complement the bounds in
\cite{costaceque2025stein} and use a synchronous de Haan--LePage coupling.
If \(\alpha>1\) and \(1\le p<\alpha\), we show that
\[
  W_{p,\|\cdot\|_\infty}(F_1,F_2)
  \le
  \Gamma\left(1-\frac p\alpha\right)^{1/p}
  \mathfrak W_{\alpha,\infty}(F_1,F_2),
\]
where
\[
  \mathfrak W_{\alpha,\infty}(F_1,F_2)
  :=
  \inf_{Z^{(1)},Z^{(2)}}
  W_{\alpha,\|\cdot\|_\infty}
  \bigl(\Law(Z^{(1)}),\Law(Z^{(2)})\bigr).
\]
We also give an \(\ell_1\)-version based on coordinatewise transport costs and
an angular version in terms of canonical angular representers. These
Wasserstein bounds are especially useful when the max-stable laws have
finite \(p\)-moments, that is, when \(p<\alpha\).

A third route uses the inf--argmax decomposition
\[
  V_k(x)
  =
  \sum_{j=1}^d
  \frac{1}{x_j^\alpha}\Psi_j^{(k)}(x),
  \qquad x\in(0,\infty)^d.
\]
For many parametric families the functions \(\Psi_j^{(k)}\) are explicit;
see, for example, \cite{Htilt,FalkD}.
Consequently, the bound
\[
  \Delta_V
  =
  \sup_{\substack{x\in(0,\infty)^d\\ \min_j x_j=1}}
  \left|
    \sum_{j=1}^d
    x_j^{-\alpha}
    \bigl(
      \Psi_j^{(1)}(x)-\Psi_j^{(2)}(x)
    \bigr)
  \right|
\]
yields direct Kolmogorov estimates in terms of model parameters. We apply
this idea in particular to Brown--Resnick and Hüsler--Reiss distributions.

We also study the case where the extreme-value copula is fixed and only the
Fréchet index of the margins varies. In that setting we obtain two-sided
Kolmogorov bounds depending only on the one-dimensional Fréchet comparison,
as well as the exact formula
\[
  W_{1,\|\cdot\|_1}
  \bigl(F_{\alpha_1},F_{\alpha_2}\bigr)
  =
  d\int_0^\infty
  \left|
    t^{-1/\alpha_1}-t^{-1/\alpha_2}
  \right|e^{-t}\,dt,
  \qquad \alpha_1,\alpha_2>1.
\]
Thus, in the fixed-copula case, the Wasserstein distance is controlled
entirely by the marginal transformation.

Finally, we discuss applications to Archimax and clustered Archimax copulas
\cite{charpentier2014multivariate,chatelain2025clustered} and to
Brown--Resnick/Hüsler--Reiss models.

\subsection*{Organisation of the paper}

Section~\ref{sec:prelim} introduces the de Haan, angular, gauge-sphere,
\(\Psi\)-, and copula representations used throughout. Section~\ref{sec:main}
contains the Kolmogorov and Wasserstein bounds and the fixed-copula
comparison with varying Fréchet index. Section~\ref{sec:applications}
specialises these results to Archimax and clustered Archimax copulas and to
Brown--Resnick/Hüsler--Reiss models.

\section{Preliminaries}
\label{sec:prelim}

\paragraph{\emph{Notation.}}
Throughout, \(d\ge1\), \(\TT=\{1,\ldots,d\}\), and \(\alpha>0\).
For \(0<\alpha<\infty\), write
\[
  \norm{x}_\alpha
  :=
  \left(\sum_{j=1}^d |x_j|^\alpha\right)^{1/\alpha},
  \qquad x\in\R^d,
\]
and
\[
  \norm{x}_\infty:=\max_{1\le j\le d}|x_j|.
\]
For \(0<\alpha<1\), \(\norm{\cdot}_\alpha\) is understood as a gauge rather
than a norm. A generic norm used for polar coordinates is denoted by
\(\norm{\cdot}\), with positive unit sphere
\[
  \Sgen:=\{s\in[0,\infty)^d:\norm{s}=1\}.
\]
\subsection*{Distances}
Let \((\mathsf X,\rho)\) be a Polish metric space. For probability measures
\(\mu,\nu\) with finite \(p\)-th \(\rho\)-moment, \(p\ge1\), define
\[
  W_{p,\rho}(\mu,\nu)
  :=
  \inf_{\pi\in\Pi(\mu,\nu)}
  \left(
    \int_{\mathsf X\times\mathsf X}
    \rho(x,y)^p\,\pi(dx,dy)
  \right)^{1/p}.
\]
Here \(\Pi(\mu,\nu)\) denotes the set of couplings of \(\mu\) and \(\nu\).
We write \(\mathcal P_p(\mathsf X)\) for the probability measures with finite
\(p\)-th \(\rho\)-moment. For \(p=1\), the Kantorovich--Rubinstein duality
gives
\[
  W_{1,\rho}(\mu,\nu)
  =
  \sup_{\mathrm{Lip}_\rho(f)\le1}
  \left|
    \int f\,d\mu-\int f\,d\nu
  \right|.
\]
For random elements \(U,V\), write
\[
  W_{p,\rho}(U,V)
  :=
  W_{p,\rho}\bigl(\Law(U),\Law(V)\bigr).
\]
See, for example, \cite[Chapter~6]{villani2009}.

On \(\R_+^d\) we use the abbreviation
\[
  \dW(\mu,\nu)
  :=
  W_{1,\|\cdot\|_\infty}(\mu,\nu),
  \qquad
  \dW(U,V)
  :=
  \dW\bigl(\Law(U),\Law(V)\bigr).
\]
For \(\mu,\nu\in\mathcal P_2(\R^d)\), let
\[
  \dWd_2(\mu,\nu)
  :=
  W_{2,\|\cdot\|_2}(\mu,\nu).
\]
Then
\[
  \dW(\mu,\nu)
  \le
  W_{1,\|\cdot\|_2}(\mu,\nu)
  \le
  \dWd_2(\mu,\nu),
\]
by \(\|x\|_\infty\le\|x\|_2\) and Jensen's inequality.

For finite Borel measures \(\mu,\nu\) on the same measurable space, we use
the convention
\[
  \|\mu-\nu\|_{\mathrm{TV}}
  :=
  \sup_A|\mu(A)-\nu(A)|,
\]
where the supremum is over measurable sets \(A\).
This is the setwise convention. When \(\mu\) and \(\nu\) have the same total
mass, it equals one half of the full variation norm
\(|\mu-\nu|(\mathsf X)\). This normalisation is used in every comparison
below.

\subsection*{Max-stable laws and angular measures}

Let \(Z=(Z_1,\ldots,Z_d)\) be a nonnegative random vector satisfying
\[
  \E{Z_j^\alpha}=1,
  \qquad j=1,\ldots,d,
\]
and
\[
  \pk*{\max_{1\le j\le d}Z_j>0}=1.
\]
The associated max-stable exponent function is
\[
  V(x)
  :=
  \E*{
    \max_{1\le j\le d}
    \frac{Z_j^\alpha}{x_j^\alpha}
  },
  \qquad x\in(0,\infty)^d,
\]
and
\[
  F(x):=\exp\{-V(x)\}
\]
has unit-\(\alpha\)-Fréchet margins. Equivalently, \(F\) admits the de Haan
representation
\[
  X_j
  =
  \max_{i\ge1}
  \frac{\overline Z^{(i)}_j}{\Gamma_i^{1/\alpha}},
  \qquad j=1,\ldots,d,
\]
where \((\overline Z^{(i)})_{i\ge1}\) are i.i.d.\ copies of \(Z\), and
\(\Gamma_i=E_1+\cdots+E_i\) with \(E_i\) i.i.d.\ standard exponential
variables independent of the marks.
This spectral representation goes back to de Haan; see
\cite{deHaan,Haan2006}.

The corresponding tail measure \(\mu\) is
\[
  \mu(A)
  :=
  \int_0^\infty
  \pk{rZ\in A}\,
  \alpha r^{-\alpha-1}\,dr.
\]
It is homogeneous of order \(-\alpha\):
\[
  \mu(cA)=c^{-\alpha}\mu(A),
  \qquad c>0.
\]
Using the polar decomposition \(z=rs\), with
\[
  r=\norm{z},
  \qquad
  s=z/\norm{z}\in\Sgen,
\]
the tail measure factorises as
\[
  \mu(dr,ds)
  =
  \alpha r^{-\alpha-1}\,dr\,H(ds),
\]
where the angular measure \(H\) on \(\Sgen\) is
\[
  H(B)
  =
  \E*{
    \norm{Z}^{\alpha}
    \ind*{Z/\norm{Z}\in B}
  },
  \qquad B\in\mathcal B(\Sgen).
\]
It satisfies the marginal constraints
\[
  \int_{\Sgen}s_j^\alpha\,H(ds)=1,
  \qquad j=1,\ldots,d,
\]
and
\[
  V(x)
  =
  \int_{\Sgen}
  \max_{1\le j\le d}
  \left(\frac{s_j}{x_j}\right)^\alpha
  H(ds).
\]
For background on tail measures, angular measures, and their polar
factorisations, see \cite{stoev2005extremal,ressel2013homogeneous}.

Let
\[
  b:=H(\Sgen)=\E{\norm{Z}^{\alpha}}.
\]
If \(\Theta\sim H/b\), then
\[
  Z^*:=b^{1/\alpha}\Theta
\]
is again a de Haan representer of \(F\). We call \(Z^*\) the canonical
representer associated with the chosen polar decomposition. The max-stable
law and the tail measure are intrinsic, whereas \(H\) and \(Z^*\) depend on
the choice of sphere.

\subsection*{Positive homogeneous gauge spheres}
\label{subsec:positive-gauge-spheres}

The polar construction also works for homogeneous gauges that are not norms.
Call a Borel-measurable, \(1\)-homogeneous map
\(\tau:[0,\infty)^d\to[0,\infty)\) \emph{admissible} if there are constants
\(0<a\le b<\infty\) such that
\[
  a\|x\|_\infty
  \le
  \tau(x)
  \le
  b\|x\|_\infty,
  \qquad x\in[0,\infty)^d.
\]
Set
\[
  S_\tau:=\{s\in[0,\infty)^d:\tau(s)=1\}.
\]
For a de Haan representer \(Z\), define its \(\tau\)-angular measure by
\[
  H_\tau(B)
  :=
  \E*{
    \tau(Z)^\alpha
    \ind*{Z/\tau(Z)\in B}
  },
  \qquad B\in\mathcal B(S_\tau).
\]
It is finite because
\[
  H_\tau(S_\tau)
  =
  \E{\tau(Z)^\alpha}
  \le
  b^\alpha\sum_{j=1}^d\E{Z_j^\alpha}
  =
  b^\alpha d.
\]
Moreover,
\[
  \int_{S_\tau}s_j^\alpha\,H_\tau(ds)=1,
  \qquad j=1,\ldots,d,
\]
and
\[
  V(x)
  =
  \int_{S_\tau}
  \max_{1\le j\le d}
  \left(\frac{s_j}{x_j}\right)^\alpha
  H_\tau(ds).
\]
Indeed, both identities follow directly by substituting
\(s=Z/\tau(Z)\). If \(b_\tau:=H_\tau(S_\tau)\) and
\(\Theta_\tau\sim H_\tau/b_\tau\), then
\[
  Z^\tau:=b_\tau^{1/\alpha}\Theta_\tau
\]
is the canonical de Haan representer associated with this gauge sphere.

The main example is
\[
  \tau_\alpha(s)
  :=
  \left(\sum_{j=1}^d s_j^\alpha\right)^{1/\alpha},
  \qquad \alpha>0.
\]
It is an admissible gauge even when \(0<\alpha<1\), and
\[
  S_{\tau_\alpha}
  =
  \left\{
    s\in[0,\infty)^d:
    \sum_{j=1}^d s_j^\alpha=1
  \right\}.
\]
On this sphere,
\[
  \sup_{s\in S_{\tau_\alpha}}\|s\|_\infty^\alpha=1.
\]
This gauge formulation is a direct instance of polar disintegration for
homogeneous tail measures; see \cite{ressel2013homogeneous}.

The stable tail dependence function associated with \(F\) is
\[
  \ell(z)
  :=
  V(z^{-1/\alpha}),
  \qquad z\in[0,\infty)^d,
\]
with the convention \(0^{-1/\alpha}=+\infty\). Equivalently,
\[
  \ell(z)
  =
  \E*{
    \max_{1\le j\le d}
    z_j Z_j^\alpha
  }
  =
  \int_{\Sgen}
  \max_{1\le j\le d}
  z_j s_j^\alpha
  H(ds).
\]
Stable tail dependence functions and their homogeneous spectral
representations are discussed, for example, in
\cite{ressel2013homogeneous,FalkD}.

\subsection*{The \texorpdfstring{\(\Psi\)}{Psi}-decomposition}

We shall use the following decomposition of the exponent function. It is a
maximiser-weight version of the inf--argmax representation in
\cite[Eq.~(6.10)]{Htilt}.

\begin{lemma}
\label{lem:Psi-general}
Let
\[
  \Delta_{d-1}
  :=
  \left\{
    w\in[0,1]^d:
    \sum_{i=1}^d w_i=1
  \right\}.
\]
Let
\[
  W=(W_1,\ldots,W_d):
  (0,\infty)^d\times[0,\infty)^d
  \longrightarrow
  \Delta_{d-1}
\]
be a Borel-measurable maximiser-weight map satisfying, for every
\(x\in(0,\infty)^d\) and \(z\in[0,\infty)^d\),
\[
  W_i(x,z)>0
  \quad\Longrightarrow\quad
  i\in\arg\max_{1\le j\le d}\frac{z_j}{x_j},
\]
and
\[
  W_i(cx,z)=W_i(x,z),
  \qquad c>0.
\]
Define
\[
  \Psi_i(x)
  :=
  \E*{Z_i^\alpha W_i(x,Z)},
  \qquad i=1,\ldots,d.
\]
Then
\[
  V(x)
  =
  \sum_{i=1}^d
  \frac{1}{x_i^\alpha}\Psi_i(x),
  \qquad x\in(0,\infty)^d.
\]
Moreover,
\[
  0\le \Psi_i(x)\le1,
  \qquad
  \Psi_i(cx)=\Psi_i(x),
  \qquad c>0.
\]
\end{lemma}

\begin{proof}
For fixed \(x\) and \(z\),
\[
  \sum_{i=1}^d
  \frac{z_i^\alpha}{x_i^\alpha}W_i(x,z)
  =
  \max_{1\le j\le d}
  \frac{z_j^\alpha}{x_j^\alpha},
\]
because \(W_i(x,z)\) charges only maximisers and sums to one. Taking
expectations gives the decomposition. The bounds and the homogeneity follow
immediately from the definition of \(W\) and the normalisation
\(\E{Z_i^\alpha}=1\).
\end{proof}

The deterministic inf--argmax rule gives the Borel maximiser-weight map
\[
  W_i(x,z)
  =
  \ind*{
    i=\inf\arg\max_{1\le j\le d}\frac{z_j}{x_j}
  }.
\]
For this selector, define
\[
  C_i(x)
  :=
  \left\{
    s\in\Sgen:
    \frac{s_i}{x_i}>\frac{s_j}{x_j}\ \forall j<i,
    \quad
    \frac{s_i}{x_i}\ge\frac{s_j}{x_j}\ \forall j>i
  \right\}.
\]
Then
\[
  \Psi_i(x)
  =
  \int_{C_i(x)}s_i^\alpha\,H(ds).
\]
If ties at the maximum have zero probability under the law of \(Z\) for a
given \(x\), then
\[
  \Psi_i(x)
  =
  \E*{
    Z_i^\alpha
    \ind*{
      \frac{Z_i}{x_i}>
      \frac{Z_j}{x_j}
      \ \forall j\ne i
    }
  }.
\]
Equivalently, if \(\widetilde Z^{(i)}\) has the \(Z_i^\alpha\)-tilted law of
\(Z\), namely
\[
  \pk*{\widetilde Z^{(i)}\in A}
  =
  \E*{Z_i^\alpha\ind{Z\in A}},
\]
then, in the no-tie case,
\[
  \Psi_i(x)
  =
  \pk*{
    \frac{\widetilde Z_i^{(i)}}{x_i}>
    \frac{\widetilde Z_j^{(i)}}{x_j}
    \ \forall j\ne i
  }.
\]
When ties have positive probability, the individual functions \(\Psi_i\)
may depend on the selector, but the sum
\[
  \sum_{i=1}^d x_i^{-\alpha}\Psi_i(x)=V(x)
\]
does not.

\subsection*{Extreme-value copulas}

Let \(F\) be max-stable with unit-\(\alpha\)-Fréchet margins and stable tail
dependence function \(\ell\). Its extreme-value copula is
\[
  C(u)
  =
  \exp\{-\ell(-\log u)\},
  \qquad u\in(0,1]^d,
\]
where \(-\log u\) is understood componentwise. In the unit-Fréchet case
\(\alpha=1\), for \(u\in(0,1)^d\), this can also be written as
\[
  C(u)
  =
  \exp\{-V(-1/\log u)\}.
\]
The latter expression extends to \(u\in(0,1]^d\) by continuity.
See \cite{Haan2006,FalkD} for background on extreme-value copulas and stable
tail dependence functions.

If
\[
  V(x)
  =
  \sum_{i=1}^d x_i^{-1}\Psi_i(x),
  \qquad x\in(0,\infty)^d,
\]
is a \(\Psi\)-decomposition in the unit-Fréchet case, then, for
\(u\in(0,1)^d\),
\[
  C(u)
  =
  \prod_{i=1}^d
  u_i^{\Psi_i(-1/\log u)}.
\]
Here \(-1/\log u\) is understood componentwise. The identity extends to
\(u\in(0,1]^d\) by continuity, with boundary values interpreted as limits
from the open cube.

\section{Main results}
\label{sec:main}

Throughout this section, unless stated otherwise, \(F_k=e^{-V_k}\),
\(k=1,2\), are \(d\)-dimensional max-stable distribution functions with
unit-\(\alpha\)-Fréchet margins, \(\alpha>0\). Thus
\[
  F_{k,j}(x)=\exp\{-x^{-\alpha}\},
  \qquad x>0,\quad j=1,\ldots,d.
\]
Let \(S\) be either the positive unit sphere of a fixed norm or an admissible
gauge sphere as in Subsection~\ref{subsec:positive-gauge-spheres}, and let
\(H_k\) denote the angular measure of \(F_k\) on this common sphere. We
assume the normalisation
\[
  \int_S s_j^\alpha\,H_k(ds)=1,
  \qquad j=1,\ldots,d,\quad k=1,2.
\]
Equivalently, if \(Z^{(k)}\) is a de Haan representer of \(F_k\), then
\[
  \E{(Z_j^{(k)})^\alpha}=1,
  \qquad j=1,\ldots,d,
\]
and
\[
  V_k(x)
  =
  \E*{
    \max_{1\le j\le d}
    \frac{(Z_j^{(k)})^\alpha}{x_j^\alpha}
  }
  =
  \int_S
  \max_{1\le j\le d}
  \left(\frac{s_j}{x_j}\right)^\alpha
  H_k(ds),
  \qquad x\in(0,\infty)^d.
\]

For \(t\ge0\), define
\begin{equation}\label{eq:omega-def}
  \omega(t)
  :=
  \begin{cases}
  0, & t=0,\\[1mm]
  \displaystyle
  \frac{t}{1+t}(1+t)^{-1/t}, & t>0.
  \end{cases}
\end{equation}
Equivalently, for \(t>0\),
\[
  \omega(t)
  =
  \sup_{q>0}\{e^{-q}-e^{-(1+t)q}\}.
\]
Indeed, the supremum is attained at \(q=\log(1+t)/t\). The function
\(\omega\) is increasing on \([0,\infty)\), and we shall repeatedly use the
elementary estimate
\begin{equation}\label{eq:omega-linear}
  \omega(t)\le\frac{t}{e},
  \qquad t\ge0.
\end{equation}
For \(t>0\), this follows from
\[
  e^{-q}-e^{-(1+t)q}
  \le tqe^{-q}
  \le \frac{t}{e},
  \qquad q>0.
\]

\subsection{Kolmogorov bounds}
\label{subsec:main-K}

Our first result bounds the Kolmogorov distance by either a Wasserstein
distance between de Haan representers or a total-variation distance between
angular measures.
Admissible gauges and their angular measures were defined in
Subsection~\ref{subsec:positive-gauge-spheres}.

\begin{theorem}[Kolmogorov bounds via representers and angular measures]
\label{prop:epsfree-alpha-strong}
Let \(F_k=e^{-V_k}\), \(k=1,2\), be \(d\)-dimensional max-stable
distribution functions with unit-\(\alpha\)-Fréchet margins, \(\alpha>0\).
Let \(H_k\) be their angular measures on a common positive sphere \(S\),
normalised by
\[
  \int_S s_j^\alpha\,H_k(ds)=1,
  \qquad j=1,\ldots,d,\quad k=1,2.
\]
Here \(S\) is either a positive norm sphere or \(S_\tau\) for an admissible
gauge \(\tau\); the gauge-sphere angular construction is recalled in
Subsection~\ref{subsec:positive-gauge-spheres}. Define
\[
  \mathcal W_\alpha
  :=
  \inf_{Z^{(1)},Z^{(2)}}
  W_{1,\|\cdot\|_\infty}
  \bigl((Z^{(1)})^\alpha,(Z^{(2)})^\alpha\bigr),
  \qquad
  M_\infty:=\sup_{s\in S}\|s\|_\infty^\alpha,
\]
where the infimum is over all de Haan representers of \(F_1\) and \(F_2\).
With
\[
  \|H_1-H_2\|_{\mathrm{TV}}
  :=
  \sup_{A\in\mathcal B(S)}|H_1(A)-H_2(A)|,
\]
one has
\[
  \dK(F_1,F_2)
  \le
  \omega\left(
    \mathcal W_\alpha
    \wedge
    M_\infty\|H_1-H_2\|_{\mathrm{TV}}
  \right)
  \le
  \frac1e
  \left(
    \mathcal W_\alpha
    \wedge
    M_\infty\|H_1-H_2\|_{\mathrm{TV}}
  \right).
\]
For any fixed pair of normalised de Haan representers, one also has
\[
  \dK(F_1,F_2)
  \le
  \omega\left(
    W_{1,\|\cdot\|_\infty}
    \bigl((Z^{(1)})^\alpha,(Z^{(2)})^\alpha\bigr)
  \right)
  \le
  \frac1e
  W_{1,\|\cdot\|_\infty}
  \bigl((Z^{(1)})^\alpha,(Z^{(2)})^\alpha\bigr).
\]
\end{theorem}

\begin{proof}
Since both \(F_1\) and \(F_2\) have unit-\(\alpha\)-Fréchet margins, they
vanish whenever at least one coordinate is non-positive. Hence it is enough
to take the supremum over \(x\in(0,\infty)^d\).

We first record a scalar identity. For \(a,b>0\), define
\[
  \kappa(a,b):=\sup_{q>0}|e^{-qa}-e^{-qb}|.
\]
If \(a=b\), then \(\kappa(a,b)=0\). If \(a\ne b\), write
\[
  m:=a\wedge b,
  \qquad
  M:=a\vee b.
\]
A direct maximisation of \(q\mapsto e^{-qm}-e^{-qM}\) gives
\[
  \kappa(a,b)
  =
  \frac{M-m}{M}
  \left(\frac{m}{M}\right)^{m/(M-m)}.
\]
Equivalently,
\begin{equation}\label{eq:scalar-kappa-omega}
  \kappa(a,b)
  =
  \omega\left(\frac{|a-b|}{a\wedge b}\right).
\end{equation}

Now let \(x\in(0,\infty)^d\). Write
\[
  x=ru,
  \qquad
  r:=\min_{1\le j\le d}x_j,
  \qquad
  \min_{1\le j\le d}u_j=1.
\]
By homogeneity of the exponent functions,
\[
  V_k(x)=r^{-\alpha}V_k(u),
  \qquad k=1,2.
\]
As \(r\) ranges over \((0,\infty)\), so does \(q:=r^{-\alpha}\). Therefore,
using \eqref{eq:scalar-kappa-omega},
\[
  \dK(F_1,F_2)
  =
  \sup_{\substack{u\in(0,\infty)^d\\ \min_j u_j=1}}
  \omega\left(
    \frac{|V_1(u)-V_2(u)|}{V_1(u)\wedge V_2(u)}
  \right).
\]

For such \(u\), choose \(j_0\) with \(u_{j_0}=1\). Then, by the angular
representation on \(S\) and the marginal normalisation,
\[
  V_k(u)
  =
  \int_S
  \max_{1\le j\le d}
  \left(\frac{s_j}{u_j}\right)^\alpha
  H_k(ds)
  \ge
  \int_S s_{j_0}^\alpha\,H_k(ds)
  =
  1.
\]
Hence
\[
  \frac{|V_1(u)-V_2(u)|}{V_1(u)\wedge V_2(u)}
  \le
  |V_1(u)-V_2(u)|.
\]
Since \(\omega\) is increasing, if
\[
  \Delta_V
  :=
  \sup_{\substack{u\in(0,\infty)^d\\ \min_j u_j=1}}
  |V_1(u)-V_2(u)|,
\]
then
\begin{equation}\label{eq:K-bound-by-DeltaV}
  \dK(F_1,F_2)\le \omega(\Delta_V).
\end{equation}

We now bound \(\Delta_V\) in two ways.

First, fix arbitrary de Haan representers \(Z^{(1)}\) and \(Z^{(2)}\), and
put
\[
  Y^{(k)}:=(Z^{(k)})^\alpha,
  \qquad k=1,2.
\]
For \(u\in(0,\infty)^d\) with \(\min_j u_j=1\), define
\[
  f_u(y):=\max_{1\le j\le d}\frac{y_j}{u_j^\alpha},
  \qquad y\in[0,\infty)^d.
\]
Since \(u_j^{-\alpha}\le1\) for all \(j\), the function \(f_u\) is
\(1\)-Lipschitz with respect to \(\|\cdot\|_\infty\). Therefore,
\[
  |V_1(u)-V_2(u)|
  \le
  W_{1,\|\cdot\|_\infty}
  \bigl(
    (Z^{(1)})^\alpha,
    (Z^{(2)})^\alpha
  \bigr).
\]
Taking the supremum over \(u\) gives, for this fixed pair of representers,
\[
  \Delta_V
  \le
  W_{1,\|\cdot\|_\infty}
  \bigl(
    (Z^{(1)})^\alpha,
    (Z^{(2)})^\alpha
  \bigr).
\]
Taking the infimum over all de Haan representers gives
\begin{equation}\label{eq:DeltaV-W-bound}
  \Delta_V\le \mathcal W_\alpha.
\end{equation}

Second, by the angular representation on \(S\),
\[
  V_k(u)
  =
  \int_S
  g_u(s)H_k(ds),
  \qquad
  g_u(s):=
  \max_{1\le j\le d}
  \left(\frac{s_j}{u_j}\right)^\alpha.
\]
For \(\min_j u_j=1\) and \(s\in S\),
\[
  0\le g_u(s)
  \le
  \|s\|_\infty^\alpha
  \le
  M_\infty.
\]
Let \(\delta:=H_1-H_2\). By the layer-cake representation,
\[
  g_u(s)
  =
  \int_0^{M_\infty}
  \mathbf 1_{\{g_u(s)>t\}}\,dt.
\]
Hence
\[
\begin{aligned}
  |V_1(u)-V_2(u)|
  &=
  \left|
    \int_S g_u(s)\,\delta(ds)
  \right| \\
  &=
  \left|
    \int_0^{M_\infty}
    \delta\bigl(\{s\in S:g_u(s)>t\}\bigr)\,dt
  \right| \\
  &\le
  M_\infty\|H_1-H_2\|_{\mathrm{TV}}.
\end{aligned}
\]
Thus
\begin{equation}\label{eq:DeltaV-TV-bound}
  \Delta_V
  \le
  M_\infty\|H_1-H_2\|_{\mathrm{TV}}.
\end{equation}

Combining \eqref{eq:K-bound-by-DeltaV}, \eqref{eq:DeltaV-W-bound}, and
\eqref{eq:DeltaV-TV-bound} proves the first displayed bound of the theorem.
The linear bound follows from \eqref{eq:omega-linear}. The fixed-representer
bound follows from the fixed-representer estimate above before taking the
infimum over representers.
\end{proof}

Related angular-measure estimates with an explicit dimension factor are
given in \cite{costaceque2025stein}. To compare angular parametrisations,
take their reference sphere to be the positive simplex and push \(H_k\)
forward under
\[
  s\longmapsto(s_1^\alpha,\ldots,s_d^\alpha).
\]
The resulting measures have mass \(d\), satisfy their moment constraints,
and yield the same exponent functions. Although their symbol
\(d_{\mathrm{TV}}\) is not defined separately, the proof of their
angular-measure proposition uses it for the full \(L^1\)-variation. Under
that normalisation, their coefficient \(d\) becomes \(2d\) in our setwise
convention; under the half-\(L^1\) convention, it would remain \(d\). In
either convention their estimate has an explicit dimension factor. For the
unnormalised angular measures used here, the coefficient in our linear bound
is \(M_\infty/e\), and it equals \(1/e\) on the positive
\(\ell_\alpha\)-sphere.

\begin{remark}[Intrinsic form, \(\Psi\)-bound, and constants]
\label{rem:intrinsic-K-form}
\begin{enumerate}[(i)]
\item
The proof gives
\[
  \dK(F_1,F_2)
  =
  \sup_{\substack{u\in(0,\infty)^d\\ \min_j u_j=1}}
  \omega\left(
    \frac{|V_1(u)-V_2(u)|}{V_1(u)\wedge V_2(u)}
  \right)
  \le
  \omega(\Delta_V)
  \le
  \frac{\Delta_V}{e},
\]
where
\[
  \Delta_V
  :=
  \sup_{\substack{u\in(0,\infty)^d\\ \min_j u_j=1}}
  |V_1(u)-V_2(u)|.
\]

\item
If
\[
  V_k(x)
  =
  \sum_{j=1}^d
  x_j^{-\alpha}\Psi_j^{(k)}(x),
  \qquad k=1,2,
\]
then
\[
  \dK(F_1,F_2)
  \le
  \omega(\Delta_\Psi)
  \le
  \frac{\Delta_\Psi}{e},
\]
with
\[
  \Delta_\Psi
  :=
  \sup_{\substack{x\in(0,\infty)^d\\ \min_j x_j=1}}
  \left|
    \sum_{j=1}^d
    x_j^{-\alpha}
    \bigl(\Psi_j^{(1)}(x)-\Psi_j^{(2)}(x)\bigr)
  \right|.
\]
In particular,
\[
  \Delta_\Psi
  \le
  \sup_{\substack{x\in(0,\infty)^d\\ \min_j x_j=1}}
  \sum_{j=1}^d
  \left|\Psi_j^{(1)}(x)-\Psi_j^{(2)}(x)\right|.
\]

\item
In Theorem~\ref{prop:epsfree-alpha-strong},
\[
  M_\infty
  =
  \sup_{s\in S}\|s\|_\infty^\alpha.
\]
If \(S=\{s\in[0,\infty)^d:\|s\|=1\}\) is a positive norm sphere and
\[
  \|x\|\ge A\|x\|_\infty,
  \qquad x\in[0,\infty)^d,
\]
then
\[
  M_\infty\le A^{-\alpha}.
\]
In particular,
\[
  M_\infty=1
\]
for \(S=\{s\in[0,\infty)^d:\|s\|_p=1\}\), \(1\le p\le\infty\).

For the weaker constant
\[
  M_\alpha
  :=
  \sup_{s\in S}\|s\|_\alpha^\alpha,
\]
one has \(M_\infty\le M_\alpha\), and hence
\[
  \dK(F_1,F_2)
  \le
  \frac{M_\alpha}{e}
  \|H_1-H_2\|_{\mathrm{TV}}.
\]

\item
On the positive \(\ell_\alpha\)-gauge sphere
\[
  S_\alpha
  =
  \left\{
    s\in[0,\infty)^d:
    \sum_{j=1}^d s_j^\alpha=1
  \right\},
\]
one has
\[
  M_\infty
  =
  \sup_{s\in S_\alpha}\|s\|_\infty^\alpha
  =
  1,
  \qquad \alpha>0.
\]
Indeed,
\[
  \|s\|_\infty^\alpha
  =
  \max_{1\le j\le d}s_j^\alpha
  \le
  \sum_{j=1}^d s_j^\alpha
  =
  1,
\]
and equality is attained at the coordinate vectors. For \(0<\alpha<1\),
this uses the gauge-sphere interpretation from
Subsection~\ref{subsec:positive-gauge-spheres}.
\end{enumerate}
\end{remark}

\subsection{Wasserstein bounds via the de Haan--LePage coupling}
\label{subsec:main-Wasserstein}

The bound in the previous subsection compares exponent functions and
therefore naturally uses \(W_1\)-transport of the powered representers
\((Z^{(k)})^\alpha\). We now give Wasserstein bounds between the max-stable
laws themselves. These estimates use a different mechanism: a synchronous
de Haan--LePage coupling, which transports the raw de Haan representers
\(Z^{(k)}\) in \(W_\alpha\).

For \(p\ge1\), write \(W_{p,\|\cdot\|}\) for the \(p\)-Wasserstein distance
associated with the norm \(\|\cdot\|\). When \(F_1,F_2\) are distribution
functions, \(W_{p,\|\cdot\|}(F_1,F_2)\) denotes the Wasserstein distance
between the corresponding probability laws.

\begin{theorem}[Wasserstein bounds via the de Haan--LePage coupling]
\label{thm:Wp-maxstable-LePage}
Let \(F_k\), \(k=1,2\), be \(d\)-dimensional max-stable distribution
functions with unit-\(\alpha\)-Fréchet margins, and let \(1\le p<\alpha\).
Define
\[
  \mathfrak W_{\alpha,\infty}(F_1,F_2)
  :=
  \inf_{Z^{(1)},Z^{(2)}}
  W_{\alpha,\|\cdot\|_\infty}
  \bigl(
    \Law(Z^{(1)}),
    \Law(Z^{(2)})
  \bigr),
\]
where the infimum is over all normalised de Haan representers of \(F_1\) and
\(F_2\). Then
\begin{equation}\label{eq:Wp-infty-LePage}
  W_{p,\|\cdot\|_\infty}(F_1,F_2)
  \le
  \Gamma\left(1-\frac p\alpha\right)^{1/p}
  \mathfrak W_{\alpha,\infty}(F_1,F_2),
\end{equation}
where normalised means
\[
  \E{(Z_j^{(k)})^\alpha}=1,
  \qquad j=1,\ldots,d,\quad k=1,2.
\]

Moreover,
\begin{equation}\label{eq:Wp-l1-LePage}
  W_{p,\|\cdot\|_1}(F_1,F_2)
  \le
  \Gamma\left(1-\frac p\alpha\right)^{1/p}
  \mathfrak C_{\alpha,1}(F_1,F_2),
\end{equation}
where the coordinatewise transport functional is
\[
  \mathfrak C_{\alpha,1}(F_1,F_2)
  :=
  \inf_{Z^{(1)},Z^{(2)}}\;
  \inf_{\pi\in\Pi(\Law(Z^{(1)}),\Law(Z^{(2)}))}
  \sum_{j=1}^d
  \left(
    \int |z_j-z_j'|^\alpha\,\pi(dz,dz')
  \right)^{1/\alpha},
\]
with the outer infimum again taken over all normalised de Haan representers
of \(F_1\) and \(F_2\).
\end{theorem}

\begin{proof}
For any normalised de Haan representer \(Z\),
\[
  \E{\|Z\|_\infty^\alpha}
  \le
  \sum_{j=1}^d \E{Z_j^\alpha}
  =
  d.
\]
Hence the transport quantities appearing in the theorem are finite for every
fixed pair of representers.

Fix normalised de Haan representers \(Z^{(1)}\) and \(Z^{(2)}\) of \(F_1\)
and \(F_2\), and let \(\pi\) be a coupling of their laws. Let
\[
  (Z_i^{(1)},Z_i^{(2)})_{i\ge1}
\]
  be i.i.d.\ with law \(\pi\), and let
\[
  \Gamma_i=E_1+\cdots+E_i,
  \qquad i\ge1,
\]
where \(E_i\) are i.i.d.\ standard exponential random variables, independent
of the marks. Define
\[
  X^{(k)}
  :=
  \bigvee_{i\ge1}\Gamma_i^{-1/\alpha}Z_i^{(k)},
  \qquad k=1,2,
\]
where the maximum is componentwise. By the de Haan--LePage representation,
\(X^{(k)}\sim F_k\).

Set
\[
  D_i:=\|Z_i^{(1)}-Z_i^{(2)}\|_\infty,
\]
and let \(D\) denote a random variable with their common law.
Then
\[
  \|X^{(1)}-X^{(2)}\|_\infty
  \le
  \sup_{i\ge1}\Gamma_i^{-1/\alpha}D_i.
\]
The marked points \((\Gamma_i,D_i)\) form a Poisson point process on
\((0,\infty)\times[0,\infty)\) with intensity measure
\[
  d\gamma\,\Law(D)(dv).
\]
Therefore, for \(x>0\),
\[
\begin{aligned}
  \pk*{
    \sup_{i\ge1}\Gamma_i^{-1/\alpha}D_i\le x
  }
  &=
  \exp\left\{
    -\int_0^\infty
      \pk{\gamma^{-1/\alpha}D>x}\,d\gamma
  \right\} \\
  &=
  \exp\left\{
    -x^{-\alpha}\E{D^\alpha}
  \right\}.
\end{aligned}
\]
Thus, for \(1\le p<\alpha\),
\[
  \E*{
    \left(
      \sup_{i\ge1}\Gamma_i^{-1/\alpha}D_i
    \right)^p
  }
  =
  \Gamma\left(1-\frac p\alpha\right)
  \left(\E{D^\alpha}\right)^{p/\alpha}.
\]
Consequently,
\[
  W_{p,\|\cdot\|_\infty}(F_1,F_2)
  \le
  \Gamma\left(1-\frac p\alpha\right)^{1/p}
  \left(
    \int \|z-z'\|_\infty^\alpha\,\pi(dz,dz')
  \right)^{1/\alpha}.
\]
Taking the infimum first over \(\pi\) and then over all normalised de Haan
representers, gives \eqref{eq:Wp-infty-LePage}.

For the \(\ell_1\)-bound, use
\[
  \|X^{(1)}-X^{(2)}\|_1
  \le
  \sum_{j=1}^d
  \sup_{i\ge1}
  \Gamma_i^{-1/\alpha}
  |Z_{ij}^{(1)}-Z_{ij}^{(2)}|.
\]
By Minkowski's inequality and the same Fréchet moment computation applied
coordinatewise,
\[
\begin{aligned}
  \left(
    \E{\|X^{(1)}-X^{(2)}\|_1^p}
  \right)^{1/p}
  &\le
  \sum_{j=1}^d
  \left(
    \E*{
      \left[
        \sup_{i\ge1}
        \Gamma_i^{-1/\alpha}
        |Z_{ij}^{(1)}-Z_{ij}^{(2)}|
      \right]^p
    }
  \right)^{1/p} \\
  &=
  \Gamma\left(1-\frac p\alpha\right)^{1/p}
  \sum_{j=1}^d
  \left(
    \int |z_j-z_j'|^\alpha\,\pi(dz,dz')
  \right)^{1/\alpha}.
\end{aligned}
\]
Taking the infimum over the coupling \(\pi\) and then over all normalised
de Haan representers gives \eqref{eq:Wp-l1-LePage}.
\end{proof}

\begin{remark}[Canonical angular version]
\label{cor:canonical-angular-Wp}
Assume \(\alpha>1\), and let
\[
  S_\alpha
  :=
  \left\{
    s\in[0,\infty)^d:
    \sum_{j=1}^d s_j^\alpha=1
  \right\}.
\]
Suppose that \(H_k\), \(k=1,2\), are angular measures on \(S_\alpha\),
normalised by
\[
  \int_{S_\alpha}s_j^\alpha\,H_k(ds)=1,
  \qquad j=1,\ldots,d.
\]
Then
\[
  H_k(S_\alpha)
  =
  \int_{S_\alpha}\sum_{j=1}^d s_j^\alpha\,H_k(ds)
  =
  d.
\]
Thus \(P_k:=H_k/d\) is a probability measure, and the canonical de Haan
representer is
\[
  Z^{(k)*}=d^{1/\alpha}S^{(k)},
  \qquad S^{(k)}\sim P_k.
\]
Therefore, for \(1\le p<\alpha\),
\[
  W_{p,\|\cdot\|_\infty}(F_1,F_2)
  \le
  d^{1/\alpha}
  \Gamma\left(1-\frac p\alpha\right)^{1/p}
  W_{\alpha,\|\cdot\|_\infty}(P_1,P_2).
\]
Similarly,
\[
  W_{p,\|\cdot\|_1}(F_1,F_2)
  \le
  d^{1/\alpha}
  \Gamma\left(1-\frac p\alpha\right)^{1/p}
  \inf_{\pi\in\Pi(P_1,P_2)}
  \sum_{j=1}^d
  \left(
    \int |s_j-t_j|^\alpha\,\pi(ds,dt)
  \right)^{1/\alpha}.
\]
\end{remark}

\subsection{Fixed copula and varying Fréchet index}
\label{subsec:fixed-copula-varying-alpha}

We now revisit the comparison with a common angular measure in
\cite{costaceque2025stein}. Under the powered angular parametrisation used
there, holding the angular measure fixed is equivalent to keeping the stable
tail dependence function, and hence the extreme-value copula, fixed while
varying the Fréchet index. We complement their upper estimates with a
dependence-free two-sided Kolmogorov bound and an exact formula for the
\(\ell_1\)-Wasserstein distance. Let \(H_\ell\) be a finite Borel spectral
measure on
\(\Sgen\) satisfying
\[
  \int_{\Sgen}s_j\,H_\ell(ds)=1,
  \qquad j=1,\ldots,d.
\]
Define
\[
  \ell(z)
  :=
  \int_{\Sgen}
  \max_{1\le j\le d}s_jz_j\,H_\ell(ds),
  \qquad z\in[0,\infty)^d.
\]
Here \(H_\ell\) is a spectral measure for the fixed stable tail dependence
function \(\ell\); equivalently, it is the angular measure of the associated
unit-Fréchet law on the chosen sphere. It is held fixed as \(\alpha\) varies
and should not be confused with the \(\alpha\)-angular measure of
\(F_\alpha\) when \(\alpha\ne1\).
For \(\alpha>0\), set
\[
  V_\alpha(x)
  :=
  \ell(x_1^{-\alpha},\ldots,x_d^{-\alpha})
  =
  \int_{\Sgen}
  \max_{1\le j\le d}
  \frac{s_j}{x_j^\alpha}\,H_\ell(ds),
  \qquad x\in(0,\infty)^d,
\]
and
\[
  F_\alpha(x):=\exp\{-V_\alpha(x)\}.
\]
Then \(F_\alpha\) is max-stable with unit-\(\alpha\)-Fréchet margins, and all
\(F_\alpha\) have the same extreme-value copula, determined by \(\ell\).

For \(\alpha_1,\alpha_2>0\), define
\begin{equation}\label{eq:eta-definition}
  \eta(\alpha_1,\alpha_2)
  :=
  \sup_{q>0}
  \left|
    e^{-q^{\alpha_1}}-e^{-q^{\alpha_2}}
  \right|.
\end{equation}

\begin{theorem}[Fixed copula and varying Fréchet index]
\label{thm:fixed-copula-varying-alpha}
Let \(\alpha_1,\alpha_2>0\), and let \(F_{\alpha_1}\) and
\(F_{\alpha_2}\) be defined as above. Set
\[
  \alpha_*:=\min(\alpha_1,\alpha_2).
\]

\begin{enumerate}[{\upshape (i)}]
\item
The Kolmogorov distance satisfies
\begin{equation}\label{eq:fixed-copula-K-sandwich}
  \eta(\alpha_1,\alpha_2)
  \le
  \dK(F_{\alpha_1},F_{\alpha_2})
  \le
  1\wedge d\,\eta(\alpha_1,\alpha_2).
\end{equation}
Moreover, if
\[
  C_0:=\sup_{z>0}ze^{-z}|\log z|,
\]
then
\begin{equation}\label{eq:fixed-copula-K-linear}
  \dK(F_{\alpha_1},F_{\alpha_2})
  \le
  1\wedge
  \frac{dC_0}{\alpha_*}|\alpha_1-\alpha_2|.
\end{equation}

\item
If \(\alpha_1,\alpha_2>1\), then
\begin{equation}\label{eq:fixed-copula-W1-exact}
  W_{1,\|\cdot\|_1}
  \bigl(F_{\alpha_1},F_{\alpha_2}\bigr)
  =
  d\int_0^\infty
  \left|
    t^{-1/\alpha_1}-t^{-1/\alpha_2}
  \right|e^{-t}\,dt.
\end{equation}

\item
For \(1\le p<\alpha_*\), the power coupling yields
\begin{equation}\label{eq:fixed-copula-Wp-bound}
  W_{p,\|\cdot\|_1}
  \bigl(F_{\alpha_1},F_{\alpha_2}\bigr)
  \le
  d\,
  m_p\left(\frac1{\alpha_1},\frac1{\alpha_2}\right),
\end{equation}
where
\[
  m_p(\beta_1,\beta_2)
  :=
  \left(
    \int_0^\infty
    \left|t^{-\beta_1}-t^{-\beta_2}\right|^p e^{-t}\,dt
  \right)^{1/p}.
\]
\end{enumerate}
\end{theorem}

\begin{proof}
Let \(Y\sim F^\star\), where \(F^\star\) denotes the member of the family
corresponding to \(\alpha=1\), namely
\[
  F^\star(x)
  =
  \exp\{-\ell(x_1^{-1},\ldots,x_d^{-1})\}.
\]
For every \(\alpha>0\),
\[
  Y^{1/\alpha}\sim F_\alpha,
\]
with powers taken componentwise. Indeed, for \(x\in(0,\infty)^d\),
\[
  \pk{Y^{1/\alpha}\le x}
  =
  \pk{Y\le x^\alpha}
  =
  \exp\{-\ell(x_1^{-\alpha},\ldots,x_d^{-\alpha})\}
  =
  F_\alpha(x).
\]

We first prove the Kolmogorov bounds. For fixed \(x\in(0,\infty)^d\), define
\[
  B_i(x)
  :=
  \{Y_j\le x_j^{\alpha_i},\ j=1,\ldots,d\},
  \qquad i=1,2.
\]
Then
\[
  |F_{\alpha_1}(x)-F_{\alpha_2}(x)|
  =
  |\pk{B_1(x)}-\pk{B_2(x)}|
  \le
  \pk{B_1(x)\triangle B_2(x)}.
\]
Moreover,
\[
  B_1(x)\triangle B_2(x)
  \subset
  \bigcup_{j=1}^d
  \left\{
    x_j^{\alpha_1}\wedge x_j^{\alpha_2}
    <
    Y_j
    \le
    x_j^{\alpha_1}\vee x_j^{\alpha_2}
  \right\}.
\]
Since each margin of \(Y\) is unit-Fréchet,
\[
\begin{aligned}
  |F_{\alpha_1}(x)-F_{\alpha_2}(x)|
  &\le
  \sum_{j=1}^d
  \left|
    \pk{Y_j\le x_j^{\alpha_1}}
    -
    \pk{Y_j\le x_j^{\alpha_2}}
  \right| \\
  &=
  \sum_{j=1}^d
  \left|
    e^{-x_j^{-\alpha_1}}
    -
    e^{-x_j^{-\alpha_2}}
  \right| \\
  &\le
  d\,\eta(\alpha_1,\alpha_2).
\end{aligned}
\]
Since \(\dK(F_{\alpha_1},F_{\alpha_2})\le1\), this proves the upper bound in
\eqref{eq:fixed-copula-K-sandwich}.

For the lower bound, fix \(j\in\{1,\ldots,d\}\) and \(x>0\). Let
\(x^{(m)}\) be the vector whose \(j\)-th coordinate is \(x\) and whose other
coordinates are equal to \(m\). By continuity of \(\ell\) and the marginal
normalisation \(\ell(e_j)=1\),
\[
  F_{\alpha_i}(x^{(m)})
  =
  \exp\left\{
    -\ell\left(
      m^{-\alpha_i},\ldots,m^{-\alpha_i},
      x^{-\alpha_i},
      m^{-\alpha_i},\ldots,m^{-\alpha_i}
    \right)
  \right\}
  \longrightarrow
  \exp\{-x^{-\alpha_i}\},
  \qquad m\to\infty,
\]
for \(i=1,2\), where \(x^{-\alpha_i}\) is placed in the \(j\)-th coordinate.
Therefore
\[
  \dK(F_{\alpha_1},F_{\alpha_2})
  \ge
  \sup_{x>0}
  \left|
    e^{-x^{-\alpha_1}}
    -
    e^{-x^{-\alpha_2}}
  \right|
  =
  \eta(\alpha_1,\alpha_2).
\]
This proves \eqref{eq:fixed-copula-K-sandwich}.

For the linear estimate, fix \(q>0\) and set
\[
  f_q(\alpha):=e^{-q^\alpha}.
\]
Then
\[
  |f_q'(\alpha)|
  =
  e^{-q^\alpha}q^\alpha|\log q|.
\]
Writing \(z=q^\alpha\), we have
\[
  q^\alpha|\log q|
  =
  \frac{z|\log z|}{\alpha}.
\]
Thus, for \(\alpha\ge\alpha_*\),
\[
  |f_q'(\alpha)|
  \le
  \frac{C_0}{\alpha_*}.
\]
The mean-value theorem gives
\[
  \eta(\alpha_1,\alpha_2)
  \le
  \frac{C_0}{\alpha_*}|\alpha_1-\alpha_2|,
\]
and \eqref{eq:fixed-copula-K-linear} follows from
\eqref{eq:fixed-copula-K-sandwich}.

We next prove the exact \(W_1\)-formula. The coupling
\[
  \bigl(Y^{1/\alpha_1},Y^{1/\alpha_2}\bigr)
\]
gives
\[
\begin{aligned}
  W_{1,\|\cdot\|_1}
  \bigl(F_{\alpha_1},F_{\alpha_2}\bigr)
  &\le
  \E*{
    \left\|
      Y^{1/\alpha_1}-Y^{1/\alpha_2}
    \right\|_1
  } \\
  &=
  \sum_{j=1}^d
  \E*{
    \left|
      Y_j^{1/\alpha_1}-Y_j^{1/\alpha_2}
    \right|
  }.
\end{aligned}
\]
Since \(Y_j\) is unit Fréchet, \(T_j:=1/Y_j\) is standard exponential.
Therefore
\[
  \E*{
    \left|
      Y_j^{1/\alpha_1}-Y_j^{1/\alpha_2}
    \right|
  }
  =
  \int_0^\infty
  \left|
    t^{-1/\alpha_1}-t^{-1/\alpha_2}
  \right|e^{-t}\,dt.
\]
This proves the upper bound in \eqref{eq:fixed-copula-W1-exact}.

For the reverse inequality, let \(\pi\) be any coupling of
\(F_{\alpha_1}\) and \(F_{\alpha_2}\), and let \(\pi_j\) be its \(j\)-th
coordinate marginal. Then \(\pi_j\) is a coupling of the corresponding
one-dimensional unit-\(\alpha_i\)-Fréchet laws. Hence
\[
\begin{aligned}
  \int\|x-y\|_1\,\pi(dx,dy)
  &=
  \sum_{j=1}^d
  \int |x_j-y_j|\,\pi_j(dx_j,dy_j) \\
  &\ge
  \sum_{j=1}^d
  W_1\bigl(F_{\alpha_1}^{(1)},F_{\alpha_2}^{(1)}\bigr),
\end{aligned}
\]
where \(F_\alpha^{(1)}\) denotes the one-dimensional unit-\(\alpha\)-Fréchet
law. The one-dimensional quantile formula gives
\[
  W_1\bigl(F_{\alpha_1}^{(1)},F_{\alpha_2}^{(1)}\bigr)
  =
  \int_0^\infty
  \left|
    t^{-1/\alpha_1}-t^{-1/\alpha_2}
  \right|e^{-t}\,dt.
\]
Taking the infimum over all couplings \(\pi\) proves
\eqref{eq:fixed-copula-W1-exact}.

Finally, for \(1\le p<\alpha_*\), the same power coupling and Minkowski's
inequality yield
\[
\begin{aligned}
  W_{p,\|\cdot\|_1}
  \bigl(F_{\alpha_1},F_{\alpha_2}\bigr)
  &\le
  \left(
    \E*{
      \left\|
        Y^{1/\alpha_1}-Y^{1/\alpha_2}
      \right\|_1^p
    }
  \right)^{1/p} \\
  &\le
  \sum_{j=1}^d
  \left(
    \E*{
      \left|
        Y_j^{1/\alpha_1}-Y_j^{1/\alpha_2}
      \right|^p
    }
  \right)^{1/p} \\
  &=
  d\,
  m_p\left(\frac1{\alpha_1},\frac1{\alpha_2}\right).
\end{aligned}
\]
This proves \eqref{eq:fixed-copula-Wp-bound}.
\end{proof}

\begin{remark}[Further consequences and radial comparison]
\label{rem:fixed-copula-further}
\begin{enumerate}[(i)]
\item
Let
\[
  \beta_i:=\frac1{\alpha_i},\qquad i=1,2,
  \qquad
  \beta_-:=\min(\beta_1,\beta_2),
  \qquad
  \beta_+:=\max(\beta_1,\beta_2).
\]
For \(\alpha_1,\alpha_2>1\),
\[
  W_{1,\|\cdot\|_1}
  \bigl(F_{\alpha_1},F_{\alpha_2}\bigr)
  \le
  d\,L_1(\beta_-,\beta_+)
  \frac{|\alpha_1-\alpha_2|}{\alpha_1\alpha_2},
\]
where
\[
  L_1(\beta_-,\beta_+)
  :=
  \int_0^1 t^{-\beta_+}|\log t|e^{-t}\,dt
  +
  \int_1^\infty t^{-\beta_-}\log t\,e^{-t}\,dt.
\]
Indeed,
\[
  |t^{-\beta_1}-t^{-\beta_2}|
  \le
  |\beta_1-\beta_2|
  \begin{cases}
    t^{-\beta_+}|\log t|, & 0<t\le1,\\
    t^{-\beta_-}\log t, & t\ge1.
  \end{cases}
\]

\item
The comparison in Theorem~\ref{thm:fixed-copula-varying-alpha} keeps the
stable tail dependence function, and hence the extreme-value copula, fixed.
A different comparison instead fixes a common finite measure \(G\) on
\(\Sgen\) and lets the radial powers vary:
\[
  V_{\alpha_i}(x)
  :=
  \int_{\Sgen}
  \max_{1\le j\le d}
  \left(\frac{s_j}{x_j}\right)^{\alpha_i}
  G(ds),
  \qquad i=1,2.
\]
This comparison generally does not keep the extreme-value copula fixed. The
law \(F_i=e^{-V_{\alpha_i}}\) has unit-\(\alpha_i\)-Fréchet margins if and
only if
\[
  \int_{\Sgen}s_j^{\alpha_i}\,G(ds)=1,
  \qquad j=1,\ldots,d,\quad i=1,2.
\]
Assume these normalisations hold. For
\(u\in(0,\infty)^d\) with \(\min_j u_j=1\), set
\[
  A_i(u):=V_{\alpha_i}(u).
\]
Then \(A_i(u)\ge1\), and
\[
  \dK(F_1,F_2)
  =
  \sup_{\substack{u\in(0,\infty)^d\\ \min_j u_j=1}}
  \Omega_{\alpha_1,\alpha_2}\bigl(A_1(u),A_2(u)\bigr),
\]
where
\[
  \Omega_{\alpha_1,\alpha_2}(a,b)
  :=
  \sup_{q>0}
  \left|
    e^{-a q^{\alpha_1}}
    -
    e^{-b q^{\alpha_2}}
  \right|.
\]
Furthermore,
\[
\begin{aligned}
  \dK(F_1,F_2)
  &\le
  \eta(\alpha_1,\alpha_2)
  +
  \sup_{\substack{u\in(0,\infty)^d\\ \min_j u_j=1}}
  \omega\left(
    \frac{
      |A_1(u)^{\alpha_2/\alpha_1}-A_2(u)|
    }{
      A_1(u)^{\alpha_2/\alpha_1}\wedge A_2(u)
    }
  \right),
\end{aligned}
\]
with the analogous bound obtained by interchanging
\(\alpha_1\) and \(\alpha_2\).

Indeed, for fixed \(u\), one compares
\[
  e^{-A_1(u)q^{\alpha_1}}
  \quad\text{with}\quad
  e^{-A_1(u)^{\alpha_2/\alpha_1}q^{\alpha_2}},
\]
which gives the \(\eta(\alpha_1,\alpha_2)\) term after the change of variable
\(r=A_1(u)^{1/\alpha_1}q\). One then compares
\[
  e^{-A_1(u)^{\alpha_2/\alpha_1}q^{\alpha_2}}
  \quad\text{with}\quad
  e^{-A_2(u)q^{\alpha_2}},
\]
which gives the displayed \(\omega\)-term.
\end{enumerate}
\end{remark}

\section{Applications}
\label{sec:applications}

This section records consequences of the main bounds in settings where the
comparison can be expressed in terms of model parameters. We focus on two
classes: Archimax and clustered Archimax copulas, and
Brown--Resnick/Hüsler--Reiss distributions.

\subsection{Archimax and clustered Archimax copulas}
\label{subsec:archimax}

Assume \(d\ge2\) in this subsection. Following
\cite{charpentier2014multivariate}, let
\(\psi:[0,\infty)\to(0,1]\) be a strict \(d\)-monotone Archimedean generator
that is continuously differentiable on \((0,\infty)\). Thus \(\psi\) is
strictly decreasing and satisfies
\[
  \psi(0)=1,
  \qquad
  \lim_{t\to\infty}\psi(t)=0.
\]
Here \(d\)-monotonicity means that \(\psi\) is differentiable up to order
\(d-2\),
\[
  (-1)^m\psi^{(m)}(t)\ge0,
  \qquad m=0,\ldots,d-2,
\]
and \((-1)^{d-2}\psi^{(d-2)}\) is nonincreasing and convex. Complete
monotonicity is a convenient stronger sufficient condition. Under these
assumptions, the Archimax construction below defines a \(d\)-copula for every
\(d\)-variate stable tail dependence function.
Let \(\ell_k\), \(k=1,2\), be stable tail dependence functions. The
corresponding Archimax copulas are
\[
  C_k(u)
  =
  \psi\left(
    \ell_k\bigl(
      \psi^{-1}(u_1),\ldots,\psi^{-1}(u_d)
    \bigr)
  \right),
  \qquad
  u\in(0,1]^d.
\]
This class contains Archimedean copulas, obtained from
\[
  \ell(z)=z_1+\cdots+z_d,
\]
and extreme-value copulas, obtained from \(\psi(t)=e^{-t}\), as special
cases.

For copulas, write
\[
  \dK(C_1,C_2)
  :=
  \sup_{u\in[0,1]^d}|C_1(u)-C_2(u)|.
\]
Set
\[
  K_\psi:=\sup_{t\ge0}t|\psi'(t)|.
\]
If \(K_\psi<\infty\), then
\begin{equation}\label{eq:archimax-main-bound}
  \dK(C_1,C_2)
  \le
  K_\psi
  \sup_{\substack{w\in[0,\infty)^d\\ \max_i w_i=1}}
  |\ell_1(w)-\ell_2(w)|.
\end{equation}

Indeed, put
\[
  z_i:=\psi^{-1}(u_i),
  \qquad
  r:=\max_{1\le i\le d}z_i.
\]
If \(r=0\), then \(u=(1,\ldots,1)\) and \(C_1(u)=C_2(u)=1\). Otherwise set
\[
  w_i:=\frac{z_i}{r},
  \qquad i=1,\ldots,d.
\]
Then \(\max_i w_i=1\), and by homogeneity of \(\ell_k\),
\[
  C_k(u)=\psi(r\ell_k(w)),
  \qquad k=1,2.
\]
The mean-value theorem gives
\[
  |C_1(u)-C_2(u)|
  \le
  r|\psi'(\xi)|\,|\ell_1(w)-\ell_2(w)|,
\]
where \(\xi=r\theta\) and
\[
  \theta
  \in
  [\ell_1(w)\wedge\ell_2(w),\ell_1(w)\vee\ell_2(w)].
\]
Since \(\ell_k(w)\ge\max_i w_i=1\), we have \(\theta\ge1\), and hence
\[
  r|\psi'(\xi)|
  =
  \frac{\xi}{\theta}|\psi'(\xi)|
  \le
  \xi|\psi'(\xi)|
  \le
  K_\psi.
\]
This proves \eqref{eq:archimax-main-bound}. Boundary points of
\([0,1]^d\) are obtained by continuity.

The right-hand side of \eqref{eq:archimax-main-bound} is the same
exponent-function discrepancy that appears in
Remark~\ref{rem:intrinsic-K-form}. Indeed, stable tail dependence functions
are continuous on \([0,\infty)^d\). Hence
\[
  \sup_{\substack{w\in[0,\infty)^d\\ \max_i w_i=1}}
  |\ell_1(w)-\ell_2(w)|
  =
  \sup_{\substack{w\in(0,\infty)^d\\ \max_i w_i=1}}
  |\ell_1(w)-\ell_2(w)|,
\]
because the latter set is dense in the former. On the positive section
\(\{w\in(0,\infty)^d:\max_i w_i=1\}\), the change of variables
\[
  x_i:=\frac1{w_i},
  \qquad i=1,\ldots,d,
\]
maps this section bijectively onto
\[
  \{x\in(0,\infty)^d:\min_i x_i=1\}.
\]
If \(\ell_k(z)=V_k(1/z)\), then on this positive section
\[
  |\ell_1(w)-\ell_2(w)|
  =
  |V_1(x)-V_2(x)|.
\]
Consequently, if \(\ell_k\) correspond to unit-Fréchet max-stable laws
\(F_k=e^{-V_k}\), then
\begin{equation}\label{eq:archimax-bound-main-tools}
  \dK(C_1,C_2)
  \le
  K_\psi
  \left(
    \mathcal W_1
    \wedge
    M_\infty\|H_1-H_2\|_{\mathrm{TV}}
  \right),
\end{equation}
where \(\mathcal W_1\), \(M_\infty\), and \(H_k\) are as in
Theorem~\ref{prop:epsfree-alpha-strong} with \(\alpha=1\).
Equivalently, using the \(\Psi\)-decomposition,
\begin{equation}\label{eq:archimax-Psi-bound}
  \dK(C_1,C_2)
  \le
  K_\psi
  \sup_{\substack{x\in(0,\infty)^d\\ \min_i x_i=1}}
  \left|
    \sum_{i=1}^d
    \frac{1}{x_i}
    \bigl(\Psi_i^{(1)}(x)-\Psi_i^{(2)}(x)\bigr)
  \right|.
\end{equation}

Motivated by the clustered Archimax construction in
\cite{chatelain2025clustered}, we next consider a general blockwise
composition. Let
\(\{I_g\}_{g=1}^G\) be a partition of \(\{1,\ldots,d\}\). For
\(k=1,2\), let \(R_k\) be a \(G\)-variate copula and let \(A_{k,g}\) be a
copula on the coordinate block \(I_g\). Assume that the functions
\[
  C_k(u)
  :=
  R_k\bigl(
    A_{k,1}(u_{I_1}),\ldots,A_{k,G}(u_{I_G})
  \bigr),
  \qquad u\in[0,1]^d,
\]
are copulas. Then
\begin{equation}\label{eq:clustered-archimax-bound}
  \dK(C_1,C_2)
  \le
  \dK(R_1,R_2)
  +
  \sum_{g=1}^G
  \dK(A_{1,g},A_{2,g}).
\end{equation}
Indeed, for fixed \(u\in[0,1]^d\), set
\[
  a_k(u)
  :=
  \bigl(
    A_{k,1}(u_{I_1}),\ldots,A_{k,G}(u_{I_G})
  \bigr),
  \qquad k=1,2.
\]
Using the fact that copulas are \(1\)-Lipschitz with respect to the
\(\ell_1\)-norm, we obtain
\[
\begin{aligned}
  |C_1(u)-C_2(u)|
  &\le
  |R_1(a_1(u))-R_2(a_1(u))|
  +
  |R_2(a_1(u))-R_2(a_2(u))|  \\
  &\le
  \dK(R_1,R_2)
  +
  \sum_{g=1}^G
  |A_{1,g}(u_{I_g})-A_{2,g}(u_{I_g})|  \\
  &\le
  \dK(R_1,R_2)
  +
  \sum_{g=1}^G
  \dK(A_{1,g},A_{2,g}).
\end{aligned}
\]
Taking the supremum over \(u\in[0,1]^d\) proves
\eqref{eq:clustered-archimax-bound}.

If the block copulas \(A_{k,g}\) are Archimax copulas, each block term on the
right-hand side can be bounded by
\eqref{eq:archimax-main-bound}, \eqref{eq:archimax-bound-main-tools}, or
\eqref{eq:archimax-Psi-bound}. If the outer copulas \(R_k\) are themselves
Archimax copulas, the same bounds apply to \(\dK(R_1,R_2)\).

\subsection{Brown--Resnick and Hüsler--Reiss models}
\label{subsec:BR-HR}

Let
\[
  U^{(k)}=(U^{(k)}_1,\ldots,U^{(k)}_d)
  \sim \mathcal N_d(0,\Sigma_k),
  \qquad k=1,2,
\]
and put
\[
  \sigma_{k,i}^2:=(\Sigma_k)_{ii}.
\]
Define the lognormal de Haan representers
\[
  Z_i^{(k)}
  :=
  \exp\left(
    U_i^{(k)}-\frac12\sigma_{k,i}^2
  \right),
  \qquad i=1,\ldots,d.
\]
Then
\[
  \E{Z_i^{(k)}}=1,
\]
and
\[
  V_k(x)
  =
  \E*{
    \max_{1\le i\le d}
    \frac{Z_i^{(k)}}{x_i}
  },
  \qquad x\in(0,\infty)^d,
\]
is the exponent function of a Brown--Resnick distribution; see
\cite{bro1977,kab2009}.

The standard Brown--Resnick model has unit-Fréchet margins, corresponding to
\(\alpha=1\). Hence the LePage Wasserstein bounds for the max-stable laws
themselves do not yield \(W_1\)-bounds in this case. Nevertheless, the
representer Wasserstein distance gives a Kolmogorov bound through
Theorem~\ref{prop:epsfree-alpha-strong}.

Let
\[
  a_k:=\frac12\operatorname{diag}(\Sigma_k)
  =
  \frac12(\sigma_{k,1}^2,\ldots,\sigma_{k,d}^2),
\]
and define the shifted Gaussian vectors
\[
  \widetilde U^{(k)}:=U^{(k)}-a_k.
\]
Then \(Z^{(k)}=\exp(\widetilde U^{(k)})\), componentwise, and
\[
  \widetilde U^{(k)}
  \sim
  \mathcal N_d(-a_k,\Sigma_k).
\]
Set
\[
  L(\Sigma_1,\Sigma_2)
  :=
  \left[
    \sum_{i=1}^d
    \left(
      e^{\sigma_{1,i}^2/2}
      +
      e^{\sigma_{2,i}^2/2}
    \right)^2
  \right]^{1/2}.
\]
For any coupling of \(\widetilde U^{(1)}\) and \(\widetilde U^{(2)}\),
\[
\begin{aligned}
  \E{\|Z^{(1)}-Z^{(2)}\|_\infty}
  &\le
  \sum_{i=1}^d
  \E{\left|e^{\widetilde U_i^{(1)}}-e^{\widetilde U_i^{(2)}}\right|} \\
  &\le
  \sum_{i=1}^d
  \left(
    e^{\sigma_{1,i}^2/2}
    +
    e^{\sigma_{2,i}^2/2}
  \right)
  \left(
    \E{
      \left|
        \widetilde U_i^{(1)}-\widetilde U_i^{(2)}
      \right|^2
    }
  \right)^{1/2}  \\
  &\le
  L(\Sigma_1,\Sigma_2)
  \left(
    \E{
      \|\widetilde U^{(1)}-\widetilde U^{(2)}\|_2^2
    }
  \right)^{1/2}.
\end{aligned}
\]
The second inequality follows from
\[
  |e^a-e^b|\le(e^a+e^b)|a-b|,
\]
Cauchy--Schwarz, and the identity
\[
  \E{e^{2\widetilde U_i^{(k)}}}=e^{\sigma_{k,i}^2}.
\]
Taking the infimum over all couplings yields
\[
  W_{1,\|\cdot\|_\infty}(Z^{(1)},Z^{(2)})
  \le
  L(\Sigma_1,\Sigma_2)
  \dWd_2\bigl(
    \mathcal N_d(-a_1,\Sigma_1),
    \mathcal N_d(-a_2,\Sigma_2)
  \bigr).
\]
Therefore
\begin{equation}\label{eq:BR-lognormal-K-bound-omega}
  \dK(F_1,F_2)
  \le
  \omega\left(
    L(\Sigma_1,\Sigma_2)
    \dWd_2\bigl(
      \mathcal N_d(-a_1,\Sigma_1),
      \mathcal N_d(-a_2,\Sigma_2)
    \bigr)
  \right).
\end{equation}
In particular,
\begin{equation}\label{eq:BR-lognormal-K-bound-linear}
  \dK(F_1,F_2)
  \le
  \frac{L(\Sigma_1,\Sigma_2)}{e}
  \dWd_2\bigl(
    \mathcal N_d(-a_1,\Sigma_1),
    \mathcal N_d(-a_2,\Sigma_2)
  \bigr).
\end{equation}

For Gaussian laws, the exact quadratic Wasserstein formula
\cite{gelbrich1990formula} gives
\[
\begin{aligned}
  &\dWd_2^2\bigl(
    \mathcal N_d(-a_1,\Sigma_1),
    \mathcal N_d(-a_2,\Sigma_2)
  \bigr) \\
  &\quad =
  \frac14
  \|\operatorname{diag}(\Sigma_1)-\operatorname{diag}(\Sigma_2)\|_2^2
  +
  \mathrm{tr}\left(
    \Sigma_1+\Sigma_2
    -2(\Sigma_2^{1/2}\Sigma_1\Sigma_2^{1/2})^{1/2}
  \right).
\end{aligned}
\]
Thus \eqref{eq:BR-lognormal-K-bound-omega} is an explicit
representer-based bound in terms of the chosen Gaussian covariance matrices.

This bound depends on the chosen Gaussian representation. The
Brown--Resnick law itself depends only on the variogram matrix
\cite{kab2009}
\[
  \Gamma_k=(\gamma_{ij}^{(k)})_{i,j=1}^d,
  \qquad
  \gamma_{ij}^{(k)}
  :=
  \operatorname{Var}(U_i^{(k)}-U_j^{(k)}).
\]
The Gaussian construction ensures that \(\Gamma_k\) is a valid variogram
matrix. Conversely, a symmetric matrix \(\Gamma=(\gamma_{ij})\) with zero
diagonal is a variogram matrix only if it is conditionally negative
semidefinite:
\[
  \sum_{i,j=1}^d a_i a_j\gamma_{ij}\le0
  \quad\text{whenever}\quad
  \sum_{i=1}^d a_i=0.
\]

We use the Hüsler--Reiss convention
\cite{MR980699}
\[
  \lambda_{ij}^{(k)}:=\sqrt{\gamma_{ij}^{(k)}},
  \qquad
  \lambda_{ii}^{(k)}=0.
\]
A common alternative convention is
\(\lambda_{ij}^{\mathrm{std}}=\sqrt{\gamma_{ij}}/2\); under that convention
the first term in the Gaussian argument below is
\(\lambda_{ij}^{\mathrm{std}}\), and the coefficient of the logarithm is
\(1/(2\lambda_{ij}^{\mathrm{std}})\).

For the displayed nondegenerate formula, assume
\(\lambda_{ij}^{(k)}>0\) for \(i\ne j\). Cases with zero off-diagonal
entries follow by continuity, or by identifying the corresponding completely
dependent coordinates. Let \(I_i:=\TT\setminus\{i\}\), listed in increasing
order. For \(m\ge1\), \(q\in\R^m\), and a positive semidefinite correlation
matrix \(R\), define
\[
  \Phi_m(q;R)
  :=
  \pk{Y_r\le q_r,\ r=1,\ldots,m},
  \qquad Y\sim\mathcal N_m(0,R).
\]
This definition includes singular correlation matrices; set \(\Phi_0:=1\).

The exponent function admits the \(\Psi\)-decomposition
\[
  V_k(x)
  =
  \sum_{i=1}^d \frac{1}{x_i}\Psi_i^{(k)}(x),
\]
where
\[
  \Psi_i^{(k)}(x)
  =
  \Phi_{d-1}\bigl(q^{(k,i)}(x);R^{(k,i)}\bigr).
\]
The vector \(q^{(k,i)}(x)\) and the rows and columns of \(R^{(k,i)}\) are
indexed by \(I_i\). For \(j,l\in I_i\),
\[
  q_j^{(k,i)}(x)
  =
  \frac{\lambda_{ij}^{(k)}}{2}
  +
  \frac{1}{\lambda_{ij}^{(k)}}
  \log\frac{x_j}{x_i},
\]
and
\[
  (R^{(k,i)})_{jl}
  =
  \frac{
    (\lambda_{ij}^{(k)})^2
    +
    (\lambda_{il}^{(k)})^2
    -
    (\lambda_{jl}^{(k)})^2
  }{
    2\lambda_{ij}^{(k)}\lambda_{il}^{(k)}
  }.
\]
The matrix \(R^{(k,i)}\) is the correlation matrix of the standardised
Gaussian increments relative to coordinate \(i\), and may be singular.

For completeness, the formula follows directly from exponential tilting.
For fixed \(k,i\), define
\[
  \frac{d\mathbb P_{k,i}}{d\mathbb P}:=Z_i^{(k)}.
\]
Under \(\mathbb P_{k,i}\), the vector \(U^{(k)}\) has law
\(\mathcal N_d(\Sigma_k e_i,\Sigma_k)\). Hence the centred, standardised
increments
\[
  \left(
    \frac{
      U_j^{(k)}-U_i^{(k)}
      -[(\Sigma_k)_{ji}-(\Sigma_k)_{ii}]
    }{\lambda_{ij}^{(k)}}
  \right)_{j\in I_i}
\]
are Gaussian with correlation matrix \(R^{(k,i)}\), and
\[
  \frac{Z_i^{(k)}}{x_i}\ge\frac{Z_j^{(k)}}{x_j}
  \quad\Longleftrightarrow\quad
  \frac{
    U_j^{(k)}-U_i^{(k)}
    -[(\Sigma_k)_{ji}-(\Sigma_k)_{ii}]
  }{\lambda_{ij}^{(k)}}
  \le q_j^{(k,i)}(x).
\]
Ties have probability zero under the nondegeneracy assumption, which proves
the displayed \(\Psi\)-formula. This is the finite-dimensional
exponential-tilting representation from \cite{Htilt}; its bivariate form
goes back to \cite{joe1994multivariate}. For \(d=1\), the convention
\(\Phi_0=1\) gives \(V_k(x)=x^{-1}\).
Consequently, Remark~\ref{rem:intrinsic-K-form} yields the
variogram-invariant bound
\begin{equation}\label{eq:HR-Psi-bound}
  \dK(F_1,F_2)
  \le
  \omega(\Delta_{\mathrm{HR}})
  \le
  \frac{\Delta_{\mathrm{HR}}}{e},
\end{equation}
where
\[
  \Delta_{\mathrm{HR}}
  :=
  \sup_{\substack{x\in(0,\infty)^d\\ \min_i x_i=1}}
  \left|
    \sum_{i=1}^d
    \frac{1}{x_i}
    \left[
      \Phi_{d-1}\bigl(q^{(1,i)}(x);R^{(1,i)}\bigr)
      -
      \Phi_{d-1}\bigl(q^{(2,i)}(x);R^{(2,i)}\bigr)
    \right]
  \right|.
\]
The bound \eqref{eq:HR-Psi-bound} is intrinsic to the Hüsler--Reiss
parameters and is therefore preferable when one wants a comparison that is
invariant under changes of Gaussian representation.

\bibliographystyle{ieeetr}
\bibliography{Max_stable_DistanceV21}
\end{document}